\documentclass{elsart}
\usepackage{natbib}

\usepackage{rotating}
  \usepackage{graphicx}
  \usepackage{psfrag}
  \usepackage{amsfonts}

\usepackage{url}
\usepackage{color}
\usepackage{amssymb}

\newtheorem{theorem}{Theorem}
\newtheorem{lemma}[theorem]{Lemma}
\newtheorem{Lemma*}[theorem]{Lemma}

\newtheorem{proposition}[theorem]{Proposition}

\newtheorem{corollary}[theorem]{Corollary}
\newtheorem{example}[theorem]{Example}

\newcommand{\real}{\mathbb R}
\newcommand{\jdlZ}{\mathbb Z}
\def\jdlqed{\vbox{\hrule \hbox{\vrule\hbox to
5pt{\vbox to 6pt{\vfil}\hfil}\vrule}\hrule}}

\newcommand{\Z}{{\mathbb Z}}
\newcommand{\N}{{\mathbb N}}

\newcommand{\C}{{\mathbb C}}
\newcommand{\R}{\mathbb R}

\usepackage{amssymb}
\journal{The Journal of Symbolic Computation}

\begin{document}
\setlength{\parindent}{0pt}
\setlength{\parskip}{2ex plus 0.4ex minus 0.4ex}

\begin{frontmatter}


\title{Short Rational Functions for Toric Algebra and Applications\thanksref{vigre}}
\thanks[vigre]{Research
supported by NSF Grants DMS-0073815, DMS-0070774, 
DMS-0200729, and by NSF VIGRE Grant DMS-0135345.}
\author[ucd]{J. A. De Loera},
\author[ucd]{D. Haws},
\author[ucd]{R. Hemmecke},
\author[ucd]{P. Huggins},
\author[ucb]{B. Sturmfels}, and
\author[ucd]{R. Yoshida}

\address[ucd]{University of California at Davis, One Shields Ave. Davis, CA 95616, USA}
\address[ucb]{University of California at Berkeley, Berkeley, CA 94720, USA}






\begin{abstract}
{\small
We encode the binomials belonging to the toric ideal $I_A$ associated
with an integral $d \times n$ matrix $A$ using a short sum of rational
functions as introduced by Barvinok \cite{bar,newbar}. 
Under the assumption that $d,n$ are fixed, this representation 
allows us to compute the Graver basis and the reduced Gr\"obner
basis of the ideal $I_A$, with respect to any term order, in time polynomial
in the size of the input. We also derive a polynomial time algorithm for
normal form computation which replaces in this new encoding the usual 
reductions typical of the division algorithm. We describe other applications,
such as the computation of Hilbert series of normal semigroup rings, and
we indicate further connections to integer programming and statistics.
}
\end{abstract}


\end{frontmatter}

\section{Introduction}
\label{intro} In this note we present polynomial-time algorithms
for computing with toric ideals and semigroup rings.
For background on these algebraic objects
and their interplay with polyhedral geometry see
\citep{stanley0},  \citep{sturmfels}, \citep{villarreal}.
Our results are a direct application of recent results
by Barvinok and Woods (2003) on short encodings
of rational generating functions (such as Hilbert series). 

Let $A = (a_{ij})$ be an integral $d \times n$-matrix 
and $b \in \jdlZ^d$ such that
the convex polyhedron
$\, P \,= \,
\{\, u \in \real^n \,\, :\,\, A \cdot u = b \,\,\hbox{and}\,\, u \geq 0 \,\}\,$
is bounded. Barvinok (1994) gave an algorithm for counting
the lattice points in $P$ in  polynomial time when 
$n-d$ is a constant. The input for Barvinok's algorithm
is the binary encoding of the integers $a_{ij}$ and $b_i$,
and the output is a formula for the multivariate
generating function $f(P)=\sum_{a \in P \cap \Z ^n} x^a$ where 
$x^a$ is an abbreviation of $x_1^{a_1} x_2^{a_2}\dots x_n^{a_n}$. This 
long polynomial with exponentially many monomials  is
encoded as a much shorter sum of rational functions of the form 
\begin{equation}
\label{barvinokseries}
f(P) \quad 
= \quad  
\sum_{i \in I} \pm \frac{x^{u_i}}{(1-x^{c_{1,i}})(1-x^{c_{2,i}})\dots 
(1-x^{c_{n-d,i}})}.
\end{equation}
Barvinok and Woods (2003) developed a set of powerful manipulation 
rules for using these short rational functions in Boolean
constructions on various sets of lattice points.
In this note we apply their
techniques to problems in combinatorial commutative algebra.
Our first theorem concerns the computation of the
\emph{toric ideal} $I_A$ of the matrix $A$. This
ideal is generated by all binomials $x^u - x^v$ 
such that $A u = A v$. In what follows, we encode any set of binomials
$\,x^u - x^v\,$ in $n$ variables as the formal sum of the
corresponding monomials $\,x^u y^v\,$ in $2n$ variables 
$x_1,\ldots,x_n,y_1,\ldots,y_n$.

\begin{theorem} \label{mainsec1} 
Let $A \in \jdlZ^{d \times n}$. Assuming that $n$ and $d$ are fixed,
there is a polynomial time algorithm to compute a short
rational function $G$ which represents the reduced Gr\"obner basis
of the toric ideal $I_A$ with respect to 
any given term order $\prec$.
Given $G$ and any monomial $x^a$,  the following tasks
can be performed in polynomial time:
\begin{enumerate}
\item Decide whether $x^a$ is in normal form with respect to $G$.
\item Perform one step of the division algorithm modulo $G$.
\item Compute the normal form of $x^a$ modulo 
the Gr\"obner basis $G$.
\end{enumerate}
\end{theorem}

Our research group at UC Davis has developed a computer program,
called {\tt LattE}, which efficiently counts the lattice
points in any rational polytope by computing its Barvinok 
representation (\ref{barvinokseries}). The Gr\"obner basis
and normal form algorithms of
Theorem  \ref{mainsec1} are currently being implemented in {\tt LattE}.
It is important to note that {\bf the Gr\"obner basis} $G$
which will be output by {\tt LattE} {\bf is a rational function}.
It is not the long list of binomials produced by
all other computer algebra systems.

\begin{example} \rm
 Fix $n = 4$, $d=2$ and let
$m \geq 3$ be an arbitrary integer.
Consider inputting the matrix
$ A=\left [\begin {array}{cccc} m&m-1&1&0\\\noalign{\smallskip}0&1&m-1&m
\end {array}\right ]$ and the lexicographic term order
into {\tt LattE}.  The task is to compute
the kernel $I_A$ of 
$$ k[x_1,x_2,x_3,x_4] \,\rightarrow \,
k[s,t] \, , \,\,\,\,
x_1 \mapsto s^m , \,
x_2 \mapsto s^{m-1} t , \,
x_3 \mapsto s t^{m-1} , \,
x_4 \mapsto t^m . $$
The output produced by {\tt LattE} would consist of the  rational function
$$ G(x,y) \,\,\, = \,\,\,
{x_1}\,{x_4}\,{y_2}\,{y_3} \,\, + \,\, {x_2}\,x_4^{m-1}
y_3^{m} \,\, + \,\,\,{\frac {{x_1}\,{x_3}\,{{y_2}}^{2}\left (\left (
{x_1}\,{y_2}\right )^{m-1}-\left ({x_3}\,{y_4}\right )^{m
-1}\right )}{{x_1}\,{y_2}-{x_3}\,{y_4}}}.
$$
This rational function is  a polynomial whose number
of terms is $m+1$ and hence grows
exponentially in the size of the input. Yet,
the running time for computing $G(x,y)$ 
is bounded by a polynomial in  $\,{\rm log}(m)$. 
It is an interesting exercise to perform the
tasks (1), (2) and (3) in Theorem \ref{mainsec1} for $G(x,y)$ and 
the monomial $x_1^m x_2^m x_3^m x_4^m$.
\end{example}

The proof of Theorem \ref{mainsec1} will be given in Section 
\ref{toric}. Special attention will be paid to the Projection Theorem 
\cite[Theorem 1.7]{newbar} since projection 
of short rational functions is the most difficult step to 
implement. Its practical efficiencyÄ has yet to be investigated.
 Our proof of Theorem
\ref{mainsec1} does use the Projection Theorem, but
our Proposition \ref{nonreducedGB} in Section \ref{toric} shows
that a \emph{non-reduced} Gr\"obner basis
can be computed in polynomial time
without using the Projection Theorem.

In Section \ref{homogenized} we present what we call the {\em
homogenized Barvinok algorithm}.  This algorithm was first outlined in
\citep{latte} and it was recently implemented in {\tt LattE}. Like the
original version in \citep{bar}, it runs in polynomial time when the
dimension is fixed. But it performs much better in practice (1) when
computing the Ehrhart series of polytopes with few facets but many
vertices; (2) when computing the Hilbert series of normal semigroup
rings.  We show its effectiveness by solving the 
classical counting problems for $5 \times 5$ \emph{magic squares}
(all row, column and diagonal sums are equal) and $3 \times 3 \times 3
\times 3$ \emph{magic hypercubes} (Theorem \ref{magic}).

A \emph{normal semigroup} $S$ is the intersection of the lattice
$\jdlZ^n$ with a rational convex polyhedral cone in $\real^n$.
The \emph{Hilbert series} of $S$ is the rational generating function
$\,\sum_{ a \in S} x^a$. Barvinok and Woods (2003)
showed that this Hilbert series 
can be computed as a short rational generating function.
We show that this computation can be done without 
the Projection Theorem when the semigroup is known to be normal.

\begin{theorem} \label{mainsec2} Under the hypothesis that 
the ambient dimension $n$ is fixed, 

1) the Ehrhart series  of a rational
 convex polytope given by linear inequalities can 
be computed in polynomial time. The  Projection Theorem 
is not used in the algorithm. 

2) The same applies to computing the Hilbert series of a normal
semigroup $S$.
\end{theorem}

In the final section of the paper we sketch applications of the
above algebraic theory to Integer Programming and Statistics. These
results will be explored in detailed elsewhere.

\section{Computing Toric Ideals}  \label{toric}

We assume that the reader is familiar with toric ideals and Gr\"obner
bases  as presented in \citep{coxlittleoshea,sturmfels}.
Barvinok and Woods (2003) showed:

\begin{lemma}[Theorem 3.6 in \citep{newbar}] \label{intersect} 
Let $S_1, S_2$ be finite
subsets of $\Z^n$, for $n$ fixed. Let $f(S_1,x)$ and $f(S_2,x)$ be their
 generating functions, given as short
rational functions with at most $k$ binomials in each denominator.
Then there exist a polynomial time algorithm, which, given
$f(S_i,x)$, computes 
$$ f(S_1 \cap S_2, x) \quad =  \quad \sum_{i \in I} \gamma_i \cdot \frac { x^{u_i} } {  (1-x^{v_{i1}})  \dots (1-x^{v_{is}}) }$$
with $s \leq 2k$, where the $\gamma_i$ are rational numbers, and $u_i,v_{ij}$ nonzero
integers.
\end{lemma}

We will use this \emph{Intersection Lemma} to extract special
monomials present in the expansion of a generating function. The
essential step in the intersection algorithm is the use of the {\em
Hadamard product} \cite[Definition 3.2]{newbar} and a special monomial
substitution. The Hadamard product is a bilinear operation on
rational functions (we denote it by $*$). The computation is
carried out for pairs of summands as in (\ref{barvinokseries}). 
Note that the Hadamard product $m_1 * m_2$ of two monomials
 $m_1,m_2$ is zero unless $m_1=m_2$. We present an example 
 of computing intersections.

\begin{example} \rm
Let $S_i=\{\,x \in \R: i-2 \leq x \leq i \,\}\cap \Z$
for $i=1,2$. We rewrite their rational generating functions as in the proof of 
Theorem 3.6 in \citep{newbar}:
$f(S_1, z) = \frac{z^{-1}}{(1-z)} + \frac{z}{(1-z^{-1})} =
\frac{-z^{-2}}{(1-z^{-1})} + \frac{z}{(1-z^{-1})} = g_{11} + g_{12},$
and $f(S_2, z) = \frac{1}{(1-z)} + \frac{z^2}{(1-z^{-1})} =
\frac{-z^{-1}}{(1-z^{-1})} + \frac{z^2}{(1-z^{-1})} = g_{21} +
g_{22}$.

We need to compute four Hadamard products 
between rational functions whose denominators are products of
binomials and denominators are monomials. Lemma 3.4 in \cite{newbar}
says that, for our example, these Hadamard products are essentially 
the same as computing the functions  (\ref{barvinokseries}) of the auxiliary
polyhedron $\{(\epsilon_1,\epsilon_2) | p_1+a_1\epsilon_1=p_2+a_2\epsilon_2, \,
\epsilon_i \geq 0 \}$
where $p_1,p_2$ are the exponent of numerators of $g_{ij}'s$ involved
and $a_1,a_2$ are the exponents of the binomial denominators. For example, the
Hadamard product $g_{11} * g_{22}$ corresponds to the polyhedron
$\,\{(\epsilon_1,\epsilon_2) | \epsilon_2=4+\epsilon_1, \,
\epsilon_i \geq 0 \}$. The contribution of this 
half line is $- \frac{z^{-2}}{(1-z^{-1})}$. We find
\begin{eqnarray*}
f(S_1, z) * f(S_2, z) \quad  = &  g_{11}*g_{21} +
g_{12} * g_{22} + g_{12}*g_{21} + g_{11} * g_{22} \\
= & \quad
\frac{z^{-2}}{(1-z^{-1})} + \frac{z}{(1-z^{-1})} - 
\frac{z^{-1}}{(1-z^{-1})} - \frac{z^{-2}}{(1-z^{-1})} \\
= &
\frac{z - z^{-1}}{1 - z^{-1}} \quad = \quad 1 + z \quad  = \quad f(S_1 \cap S_2, z).
\end{eqnarray*}
\end{example}

Another key subroutine introduced by Barvinok and Woods is the
following \emph {Projection Theorem}. In both Lemmas \ref{intersect}
and \ref{project}, the dimension $n$ is assumed to be fixed.

\begin{lemma}[Theorem 1.7 in \citep{newbar}] \label{project} 
Assume the dimension $n$ is a fixed constant.
Consider  a rational polytope $P \subset \real^n$ and a
linear map $T: \Z^n \rightarrow \Z^k$.
There is a polynomial time algorithm which
computes a short representation of the
generating function $\, f \bigl(T(P \cap \Z^n),x\bigr) $.
\end{lemma}

We represent a term order $\prec$ on monomials in $x_1,\ldots,x_n$ by
an integral $n \times n $-matrix $W$ as in \citep{Mora+Robbiano}.  Two
monomials satisfy $\,x^\alpha\prec x^\beta \,$ if and only if
$W\alpha$ is lexicographically smaller than $W\beta$.  In other words,
if $w_1,\ldots,w_n$ denote the rows of $W$, there is some
$j\in\{1,\ldots,n\}$ such that $w_i\alpha=w_i\beta$ for $i<j$, and
$w_j\alpha<w_j\beta$. For example, $W=I_n$ describes the lexicographic
term ordering. In what follows, we will denote by $\prec_W$ the term
ordering defined by $W$.

\begin{lemma} \label{extractmonomial} 
Let $S \subset \Z^n_+$ be finite. Suppose the polynomial $\,f(S,x) =
\sum_{\beta \in S} x^\beta \, $ is represented as a short rational
function and let $\prec_W$ be a term order. We can extract the
(unique) leading monomial of $f(S,x)$ with respect to $\prec_W$, in
polynomial time.
\end{lemma}

\noindent {\em Proof:} The term order $\prec_W$ is represented by an
integer matrix $W$. For each of the rows $w_j$ of $W$ we perform a
monomial substitution $x_i:=x'_it_j^{w_{ji}}$.  Such a monomial
substitution can be computed in polynomial time by \cite[Theorem
2.6]{newbar}. The effect is that the polynomial $f(S,x)$ gets replaced
by a polynomial in the $t$ and the $x's$. After each substitution we
determine the degree in $t$.  This is done as follows: We want to do
calculations in univariate polynomials since this is faster so we
consider the polynomial $g(t)=f(S,1,t)$, where all variables except
$t$ are set to the constant one. Clearly the degree of $g(t)$ in $t$
is the same as the degree of $f(S,x',t)$.  We create the
\emph{interval polynomial} $i_{[p,q]}(t)=\sum^q_{i=p} t^i$ which
obviously has a short rational function representation.  Compute the
Hadamard product of and $i_{[p,q]}$ with $g(t)$. This yields those
monomials whose degree in the variable $t$ lies between $p$ and
$q$. We will keep shrinking the interval $[p,q]$ until we find the
degree.  We need a bound for the degree in $t$ of $g(t)$ to start a
binary search.  A polynomial upper bound $U$ can be found via the
estimate in Theorem 3.1 of \citep{lasserre} by easy
manipulation of the numerator and denominator of the fractions in
$g(t)$. In no more than $log(U)$ steps one can determine the degree in
$t$ of $f(S,x,t)$ by using a standard binary search algorithm.

Once the degree $r$ in $t$ is known, we compute the Hadamard product
of $f(S,x,t)$ and $i_{[r,r]}$, and then compute the limit as $t$
approaches $1$.  This can be done in polynomial time using residue
techniques. The limit represents the subseries $\,H(S,x) = \sum_{
\beta \cdot w_j = r} x^{\beta} $.  Repeat the monomial and highest
degree search for the row $w_{j+1}$,$w_{j+2}$, etc.  Since $\prec_W$
is a term order, after doing this $n$ times we will have only one
single monomial left, the desired leading monomial. \jdlqed

\begin{proposition} \label{nonreducedGB}
Let $A \in \jdlZ^{d \times n}$, $W \in \Z^{n \times n}$ specifying a term
order $\prec_W$, and assume that
$d$ and $n$ are fixed.

1) There is a polynomial time algorithm to compute a short
rational function $G$ which represents a universal Gr\"obner basis of $I_A$.

2) Given the term order $\prec_W$ and a short rational function
encoding a (possibly infinite) set of binomials $\sum x^uy^v$, one
can compute in polynomial time a short rational function encoding
only those binomials such that $x^v \prec_W x^u$.

3) Suppose we are given a sum of short rational functions $f(x)$
which is identical, in a monomial expansion, to a single monomial $x^a$.
Then in polynomial time we can recover the (unique) exponent vector $a$.
\end{proposition}

\noindent {\em Proof:}
1) Denote by $w_i$ the $i$-th row of the matrix $W$
which specifies the term order.  Set
 $\,M=(d+1)(n-d)D(A)\,$ where $D(A) $ is the
largest absolute value of any $d \times d$-subdeterminant of $A$.
Using Barvinok's algorithm in \citep{bar}, we compute the 
 following generating function in $2n$ variables:
$$ G(x,y) \quad = \quad  \sum \bigl\{
\, x^u y^v \, \, : \,\, A u = A v\, \,\,\hbox{and}
\,  \,  0 \leq u_i,v_i \leq M \, \bigr\}. $$
This is the sum over all lattice points in a rational polytope.
Lemma 4.1 and Theorem 4.7 in Chapter 4 of \citep{sturmfels}
imply that the toric ideal $I_A$ is generated by
the finite set of  binomials $x^u-x^v$ corresponding to the
terms $x^u y^v$ in $G(x,y)$.
Moreover, these binomials are  a  universal Gr\"obner basis of $I_A$.

2) Suppose we are given a short rational generating function
$\,G_0(x,y) \, = \, \sum x^u y^v \,$ representing a set
 of binomials $\,x^u - x^v$ in $I_A$, for instance
$G_0 = G$ in part (1). In the following steps, we will alter the series so that
a term $x^u y^v$ gets removed whenever
$u$ is not bigger than $v$ in the term order $\prec_W$.
Starting with $H_0 = G_0$, we perform Hadamard products
with short rational functions $f(S;x,y)$ for $ S \subset \Z^{2n}$.

Set $H_i = H_{i-1}* f(\{(u,v):w_iu=w_iv\})$, and $G_i = H_{i-1} *
f(\{(u,v):w_iu\geq w_iv+1\}).$ All monomials $x^uy^v\in G_j$ have the
property that $w_iu=w_iv$ for $i<j$, $w_ju>w_jv$, and thus $v\prec_W
u$. On the other hand, if $v\prec_W u$ then there is some $j$ such
that $w_iu=w_iv$ for $i<j$, $w_ju>w_jv$, and we can conclude that
$x^uy^v \in G_j$.  This proves that $G=G_1\cup G_2\cup\ldots\cup G_n$
encodes exactly those binomials in $G_0$ that are correctly ordered
with respect to $\prec_W$. We have proved our claim since all of the
above constructions can be done in polynomial time.

3) Given $f(x)$ we can compute in polynomial time the partial
   derivative $\partial f(x)/\partial x_i$. This puts the exponent of
   $x_i$ as a coefficient of the unique monomial. To compute the
   derivative can be done in polynomial time by the quotient and
   product derivative rules.  Each time we differentiate a short
   rational function of the form
   $$\frac{x^{b_i}}{(1-x^{c_{1,i}})(1-x^{c_{2,i}})\dots
   (1-x^{c_{d,i}})}$$ 
  we add polynomially many (binomial) factors to the
   numerator. The factors in the numerators should be expanded into
   monomials to have again summands in short rational canonical form
   $\frac{x^{b_i}}{(1-x^{c_{1,i}})(1-x^{c_{2,i}})\dots
   (1-x^{c_{d,i}})}$. Note that at most $2^d$ monomials
   appear (a constant) each time.  Finally, if we take the limit when
   all variables $x_i$ go to one we will get the desired exponent.
\jdlqed

\begin{example} \rm Using {\tt LattE} we compute the set of
all binomials of degree less than or equal $10000$ in the toric ideal
$I_A$  of the matrix $\,A \, = \, \left
[\begin {array}{cccc} 1&1&1&1\\\noalign{\smallskip}0&1&2&3
\end {array}\right ]$. This matrix represents the
\emph{Twisted Cubic Curve} in algebraic geometry.
We find that there are exactly $195281738790588958143425$ 
such binomials. Each binomial is encoded as a monomial 
$\,x_1^{u_1} x_2^{u_2} x_3^{u_3} x_4^{u_4}
    y_1^{v_1} y_2^{v_2} y_3^{v_3} y_4^{v_4} $. The computation
takes about $40$ seconds. The output is a 
 sum of $538$ simple rational functions of the form
a monomial divided by a product such as
$
\left (1-{\frac {x_{{3}}y_{{4}}}{x_{{1}}y_{{2}}}}\right ) \left (1-{
\frac {x_{{1}}x_{{4}}y_{{2}}}{x_{{3}}}}\right ) \left (1-x_{{1}}y_{{1
}}\right ) \left (1-x_{{1}}x_{{3}}{y_{{2}}}^{2}\right ) \left (1
-x_{{3}}y_{{3}}\right ) \left (1-x_{{2}}y_{{2}}\right ) $. \jdlqed
\end{example}

\noindent {\em Proof of Theorem \ref{mainsec1}:}
Proposition \ref{nonreducedGB} gives a Gr\"obner basis for the
toric ideal $I_A$ in polynomial time.  We now show  how to get
the reduced Gr\"obner basis.

Step 1. Compute the generating function which encodes all
 binomials in  $I_A$:
$$ f(x,y) \quad = \quad \sum \bigl\{ \, x^u y^v \, \, : \,\, A u = A v
\,\,\hbox{and} \,\,u,v \geq 0 \, \bigr\},$$ This computation is
similar to part 1 of Proposition \ref{nonreducedGB} except that there
is no upper bound $M$. Next we wish to remove from $f(x,y)$ all
incorrectly ordered binomials (i.e. those monomials $x^uy^v$ with $u
\prec_W v$ instead of the other way around).  We do this following
part 2 of Proposition \ref{nonreducedGB}. Abusing notation let us
still call $f(x,y)$ the resulting sum of rational functions. Let now
$g(x)$ be the projection of $f(x,y)$ onto the first group of
variables. Thus $g(x)$ is the sum over all non-standard monomials, and
it can be computed in polynomial time by Lemma \ref{project}.

\noindent Step 2. Write  $\frac{1}{1-x} = 
\prod\limits_{i=1}^{n}\frac{1}{1-x_i}$ for the generating function of all
$x$-monomials. We compute the following
{\sl Hadamard product} of $n$ rational functions in $x$:
$$ \biggl( \frac{1}{1-x} - x_1 \cdot g(x) \biggr) * 
\biggl( \frac{1}{1-x} - x_2 \cdot g(x) \biggr) * \cdots *
\biggl( \frac{1}{1-x} - x_n \cdot g(x) \biggr). $$
This is the generating function over those monomials all of whose 
proper factors are standard monomials modulo the toric ideal $I_A$.

\noindent Step 3.
Let $h(x,y)$ denote the ordinary product of the result of Step 2 with
$$ \frac{1}{1-y} - g(y) \quad = \quad \sum \bigl\{ \,y^v \,\,: \,\, v \,\,
\hbox{standard monomial modulo} \, I_A \,\bigr\}.$$ 
Thus  $h(x,y)$ is the sum of all monomials 
$x^uy^v$ such that $x^v$ is standard and $x^u$
is minimally non-standard.
 Compute the Hadamard product $G(x,y):=\, f(x,y) * h(x,y)$.
This is a short rational representation of a polynomial, namely,
it is the sum over all monomials $x^u y^v$ such that the binomial
$x^u - x^v$ is in the reduced Gr\"obner basis of $I_A$
with respect to $W$. 

We next prove claims 1 and 2.  Let $G(x,y)$ be the 
reduced Gr\"obner basis of $I_A$ encoded by the rational function above,
and let $M$ be the degree bound of Proposition \ref{nonreducedGB}.
Let $b(x,y)$ be the rational function representing $\,\{(u,v):0\leq u \leq a, \,
0\leq v\leq M \}$.  The Hadamard product $ \,\bar{G}(x,y) = 
G(x,y)* b(x,y)$ is computable in polynomial time and encodes exactly 
those binomials in $G$  that can reduce $x^a$. If $\bar{G} $ is empty
then $x^a$ is in normal form already, otherwise we use
 Lemma \ref{extractmonomial} to find
 an element $(u,v)\in\bar{G}$ and reduce $x^a$ to $x^{a-u+v}$.

It is worth noting that analytic calculations may now be used as part
of algebraic algorithms: Suppose again we wish to decide whether $x^a$
is in reduced normal form or not. Take $G(x,y)$ as before and compute
$F(x)=G(x,1)$. This can be done using monomial substitution
\citep{newbar}. Next compute the integral
\[  \frac{ 1 }{ (2 \pi i)^n} \int_{
\left| x_1 \right| = \epsilon_1 } \cdots \int_{ \left| x_d \right| =
\epsilon_d } \frac{ (x_1^{ - a_1} \cdots x_n^{ - a_n}) F(x) }{
\left( 1 - x_1 \right) \cdots \left( 1 - x_n  \right) }
\ d x \ . 
\] 
Here $ 0 < \epsilon_1, \dots , \epsilon_d < 1 $ are different numbers
such that we can expand all the $\frac{ 1 }{ 1 - x_k} $ into the power
series about $0$. It is possible to do a partial fraction
decomposition of the integrand into a sum of simple fractions. The
integral is a non-negative integer: it is
the number of ways that the monomial $x^a$ can be written in terms of the leading 
monomials of the Gr\"obner bases $G$.

We now present the algorithm for claim 3 in 
Theorem \ref{mainsec1}. A curious byproduct of representing 
Gr\"obner bases with short rational functions is that the reduction to
normal form need not be done by dividing several times anymore:

Step 4.  Let $f(x,y)$ and $g(x)$ as above and compute the Hadamard
product
$$ H(x,y) \quad := \quad f(x,y) * \biggl( \, \bigl( \frac{1}{1-x} \bigr) \cdot
\bigl( \frac{1}{1-y} - g(y) \bigr) \biggr). $$
This is the sum over all monomials $x^u y^v$ where $x^v$ is
the normal form of $x^u$.

Step 5. We use $H(x,y)$ as one would use a traditional Gr\"obner basis of 
the ideal $I_A$.
The normal form of a monomial $x^a$ is computed by forming the Hadamard product
$$ H(x,y) *  \frac{x^a}{1-y}.$$

Since this is strictly speaking a sum of rational functions equal to a
single monomial, applying Part 3 of Proposition \ref{nonreducedGB}
concludes the proof of Theorem \ref{mainsec1}.  \jdlqed

\section{Computing Normal Semigroup Rings} 
\label{homogenized}

We observed in \citep{latte} that a major practical bottleneck of the
original Barvinok algorithm in \citep{bar} is the fact that a polytope
may have too many vertices. Since originally one visits each vertex to
compute a rational function at each tangent cone, the result can be
costly. For example, the well-known polytope of semi-magic cubes in
the $ 4\times 4 \times 4$ case has over two million vertices, but only
64 linear inequalities describe the polytope.  In such cases we
propose a homogenization of Barvinok's algorithm working with a single
cone.

There is a second motivation for looking at the
homogenization. Barvinok and Woods \citep{newbar} proved that the
Hilbert series of semigroup rings can be computed in polynomial time.
We show that for {\em normal semigroup rings} this can be done
simpler and more directly,
without using the Projection Theorem.

Given a rational polytope $P$ in $\real^{n-1}$, we set
 $i(P,m)=\# \{ z \in \Z^{n-1} : z\in m P \}$. 
The {\em Ehrhart series} of $P$ is the generating function
$\,\sum_{m=0}^\infty i(P,m) t^m$.

\begin{algorithm}[Homogenized Barvinok  algorithm]  \ \ 
\label{homogbarvinok}

\noindent {\bf Input:} A  full-dimensional, rational convex polytope $P$ 
in $\real^{n-1}$ specified by  linear inequalities and 
linear equations. \ \
 
{\bf Output:} The Ehrhart series of $P$. \rm
\begin{enumerate}
\item Place the polytope $P$ into the hyperplane defined by $x_n = 1$ in $\real^n$.
Let $K$ be the $n$-dimensional cone over
$P$, that is, $K=cone(\{(p,1) : p \in P\})$.

\item Compute the polar cone $K^*$.  The normal vectors of the facets
of $K$ are exactly the extreme rays of $K^*$.  If the polytope $P$ has
far fewer facets then vertices, then the number of rays of the cone
$K^*$ is small.

\item Apply Barvinok's decomposition of $K^*$ into unimodular cones.
Polarize back each of these cones. It is known, e.g. Corollary 2.8
in \citep{BarviPom}, that by dualizing back we get a unimodular cone
decomposition of $K$. All these cones have the same dimension as $K$.
Retrieve a signed sum of multivariate rational functions which 
represents the series $ \,\sum_{a \in K \cap \Z^{n}} x^a$.

\item Replace the variables $x_i$ by $1$ for $i \leq n-1$
and output the resulting series in $t = x_n$. This can be done using
the methods in \citep{latte}.
\end{enumerate}
\end{algorithm}

We recall that one of the key steps in Barvinok's algorithm is 
that any cone can be decomposed as the signed sum of
(indicator functions of) unimodular cones.

\begin{theorem}[see \citep{bar}] \label{barvil} Fix the dimension $n$.
Then there exists a polynomial time algorithm which decomposes a
rational polyhedral cone $K \subset \real^n$ into unimodular cones
$K_i$ with numbers $\epsilon_i \in \{-1, 1\}$ such that $$f(K) \,\, = \,\,
\sum_{i \in I} \epsilon_i f(K_i) \mbox{, } \quad |I| < \infty .$$
\end{theorem}

Originally, Barvinok had pasted together such formulas, one for each
vertex of a polytope, using a result of Brion. The point is that this
can be avoided:

\noindent{\em Proof of Theorem \ref{mainsec2}:} We first prove
part (1). The algorithm solving the problems is Algorithm
\ref{homogbarvinok}. 
Steps 1 and 2 are polynomial when the dimension is
fixed. Step 3 follows from Theorem \ref{barvil}. We require a special
monomial substitution, possibly with some poles.
This can be done in polynomial time by \citep{newbar}.

Part (2): From the characterization of the integral closure of the
semigroup $S$ as the intersection of a pointed polyhedral cone with
the lattice $\Z^n$ is clear that Algorithm 1, with the modification
that the cone $K$ in question is given by the rays of the cone (the
generators of the monomial algebra). In fixed dimension one can
transfer from the the extreme rays representation of the cone or to
the facet representation of the cone in polynomial time.
\jdlqed

\begin{corollary}
Given a normal semigroup ring $R$ of fixed Krull dimension, there is a polynomial time
algorithm which decides whether $R$ is Gorenstein.
\end{corollary}

\noindent{\em Proof:} Let $R$ be a normal semigroup ring for a
semigroup of $\Z^n$. Hochster's theorem says that the normal semigroup
rings are Cohen-Macaulay domains \citep{stanley0}. Denote by $F(R,z)$
the generating function of the monomials of normal semigroup ring $R$
(computable in polynomial time by previous theorem). Then by 
Theorem 12.7 in \citep{stanley0}, it is enough to check that $F(R,z)=(-1)^n z^a
F(R,1/z)$ for some $a \in Z^n$ efficiently. The change of variables
can be done in polynomial time and thus get $F(R,1/z)$. To check
whether $F(R,z)/F(R,1/z)$ is a single polynomial we can compute the
monomial evaluation $z_i=1$ for $i=1 \dots d$. \jdlqed

Each pointed affine semigroup $S \subset \Z^n$ can be {\em graded}. This means that
there is a linear map $\,deg: S \rightarrow \N$ with
$deg(x)=0$ if and only if $x=0$.  Given a pointed graded affine
semigroup define $S_r$ to be the set of all degree $r$ elements, i.e.
$S_r=\{x \in S : deg(x)=r\}$.  The {\em Hilbert series} of $S$ is the
formal power series $H_S(t)=\sum_{k=0}^\infty \# (S_r)
t^r$. Algebraically, this is just the Hilbert series of the semigroup
ring $\C[S]$. It is a well-known property that $H_S$ is
represented by a rational function of the form

$$\frac{Q(t)}{(1-t^{d_1})(1-t^{d_2})\dots (1-t^{d_n})}$$

where $Q(t)$ is a polynomial of degree less than $d_1+\dots+d_n$ (see
Chapter 4 \citep{stanley}).  Several other methods had been tried to
compute the Hilbert series explicitly (see \citep{ahmed} for
references).  One of the most well-known challenges was that of
counting the $5 \times 5$ magic squares of magic sum $n$. Similarly
several authors had tried before to compute the Hilbert series of the
$3 \times 3 \times 3 \times 3$ semi-magic cubes.  It is not difficult
to see this is equivalent to determining an Ehrhart series. Using
Algorithm 1 we finally present the solution, which had been
inaccessible using Gr\"obner bases methods. For comparison, the reader
familiar with Gr\"obner bases computations should be aware that the
$5\times 5$ magic squares problem required a computation of a
Gr\"obner bases of a toric ideal of a matrix $A$ with 25 rows and over
4828 columns. Our attempts to handle this problem with {\tt CoCoA} and
{\tt Macaulay2} were unsuccessful.  We now give the numerator and then
the denominator of the rational functions computed with the software
{\tt LattE}:

\begin{theorem} 
\label{magic}

\item The generating function $\,\sum_{n \geq 0} f(n) t^n \,$ for the number $f(n)$ of $5\times
5$ magic squares of magic sum $n$ is given by the rational function
$p(t)/q(t)$ with denominator

{\footnotesize

$p(t)=
{t}^{76}+28\,{t}^{75}+639\,{t}^{74}+11050\,{t}^{73}+136266\,{t}^{72}+
1255833\,{t}^{71}+9120009\,{t}^{70}+54389347\,{t}^{69}+
274778754\,{t}^{68}+1204206107\,{t}^{67}+4663304831\,{t}^{66}+
16193751710\,{t}^{65}+51030919095\,{t}^{64}+ 
147368813970\,{t}^{63}+393197605792\,{t}^{62}+975980866856\,{t}^{61}+
2266977091533\,{t}^{60}+\newline 
4952467350549\,{t}^{59}+10220353765317\,{t}^{58}+
20000425620982\,{t}^{57}+37238997469701\,{t}^{56}+\newline
66164771134709\,{t}^{55}+112476891429452\,{t}^{54}+183365550921732\,{t}^{53}+
287269293973236\,{t}^{52}+433289919534912\,{t}^{51}+
630230390692834\,{t}^{50}+885291593024017\,{t}^{49}+
1202550133880678\,{t}^{48}+1581424159799051\,{t}^{47}+
2015395674628040\,{t}^{46}+2491275358809867\,{t}^{45}+\newline
2989255690350053\,{t}^{44}+3483898479782320\,{t}^{43}+
3946056312532923\,{t}^{42}+\newline
4345559454316341\,{t}^{41}+
4654344257066635\,{t}^{40}+4849590327731195\,{t}^{39}+\newline
4916398325176454\,{t}^{38}+4849590327731195\,{t}^{37}+
4654344257066635\,{t}^{36}+\newline
4345559454316341\,{t}^{35}+
3946056312532923\,{t}^{34}+3483898479782320\,{t}^{33}+\newline
2989255690350053\,{t}^{32}+2491275358809867\,{t}^{31}+
2015395674628040\,{t}^{30}+\newline
1581424159799051\,{t}^{29}+
1202550133880678\,{t}^{28}+885291593024017\,{t}^{27}+\newline
630230390692834\,{t}^{26}+433289919534912\,{t}^{25}+
287269293973236\,{t}^{24}+183365550921732\,{t}^{23}+
112476891429452\,{t}^{22}+66164771134709\,{t}^{21}+37238997469701
\,{t}^{20}+20000425620982\,{t}^{19}+10220353765317\,{t}^{18}+
4952467350549\,{t}^{17}+2266977091533\,{t}^{16}+975980866856\,{t}^
{15}+\newline
393197605792\,{t}^{14}+147368813970\,{t}^{13}+51030919095\,{t
}^{12}+16193751710\,{t}^{11}+4663304831\,{t}^{10}+1204206107\,{t}^
{9}+274778754\,{t}^{8}+54389347\,{t}^{7}+9120009\,{t}^{6}+1255833
\,{t}^{5}+136266\,{t}^{4}+11050\,{t}^{3}+639\,{t}^{2}+28\,t+1
$ \ \ \ \ \ and numerator

$ q(t)=
\left ({t}^{2}-1\right )^{10}\left ({t}^{2}+t+1\right )
^{7}\left ({t}^{7}-1\right )^{2}
\left ({t}^{6}+{t}^{3}+1\right )\left ({t}^{5}+{t}^{3}+{t}^{2}+t+1
\right )^{4}\left (1-t\right )^{3}\left ({t
}^{2}+1\right )^{4}
$.

}

The generating function $\,\sum_{n \geq 0} f(n) t^n \,$ for the number $f(n)$ of $3\times 3
\times 3 \times 3$ magic cubes with magic sum $n$ is given the
rational function $r(t)/s(t)$ where

{\footnotesize
$
{t}^{54}+150\,{t}^{51}+5837\,{t}^{48}+63127\,{t}^{45}+331124\,{t}^{42}
+1056374\,{t}^{39}+2326380\,{t}^{36}+3842273\,{t}^{33}+5055138\,{t}^{
30}+5512456\,{t}^{27}+5055138\,{t}^{24}+3842273\,{t}^{21}+2326380\,{t}
^{18}+1056374\,{t}^{15}+331124\,{t}^{12}+63127\,{t}^{9}+5837\,{t}^{6}+
150\,{t}^{3}+1
$ \ \ \ \ \ and

$
q(t)=\left ({t}^{3}+1\right )^{4}\left ({t}^{
12}+{t}^{9}+{t}^{6}+{t}^{3}+1\right )\left (1-{t}^{3}\right )^{9}
\left ({t}^{6}+{t}^{3}+1\right )
$.
}
\end{theorem}




\section{Applications}

As explained in Chapter 5 of the book \cite{sturmfels}, Gr\"obner
bases can be useful in the context of integer programming, serving as
universal test sets of families of integer programs, and in
statistics, where they can be used as the Markov basis moves used to
generated elements uniformly at random (e.g contingency tables
counting). Therefore the fact that we can compute Gr\"obner bases
and normal forms in polynomial time (under certain hypothesis) 
can then be used to prove the following results:

\begin{corollary}
Let $A \in \Z^{d \times n}$,
$b \in \Z^d$, $W \in \Z^{n}$. Assume that $d$ and $n$ are fixed.
There is a polynomial time algorithm to solve the integer
programming problem $min_{x \in P \cap \Z^n} Wx$ where $P(b)=\{x | Ax=b,
x\geq 0 \}$.
\end{corollary}

\noindent {\bf sketch of proof:} Make the cost vector 
$W$ into a term order by breaking ties of the order $m_1 > m_2$ if 
$Wm_1> Wm_2$. You can do this via lexicographic ordering.
From Chapter 5 of \cite{sturmfels} the integral optimum of $P$
can be obtained from the Gr\"obner basis obtained in Theorem 
\ref{mainsec1} and then computing the normal form of the monomial $x^b$ 
with respect to the Gr\"obner basis. Since both steps can be done
efficiently the corollary follows.

Another application is to the uniform sampling of lattice points
inside polyhedra of the form $P(b)= \{ x \in \R^d | Ax = b \mbox{, } x
\geq 0 \}$. Given $M$ be a finite set such that $M \subset \{ x \in
\Z^d| Ax = 0 \}$.  We define the graph $G_b$ such that its nodes are
all the lattice points inside of $P$ and there is an undirected edge
between a node $u$ and a node $v$ iff $u - v \in M$. In general this
graph may not be connected for arbitrary choices of $M$. We say  
$M$ is a {\em {Markov basis}} if $G_b$ is a connected graph 
for all $b$.

\begin{corollary}
Given $A \in \Z^{d \times n}$, where $d$ and $n$ are fixed, there is a
polynomial time algorithm to compute a multivariate rational generating
function for a Markov basis $M$ associated to $A$. This is
presented as a short sum of rational functions.
\end{corollary}

We conclude with another with numeric question. Ian Dinwoodie
communicated to us the problem of counting all $7 \times 7$
contingency tables whose entries are nonnegative integers $x_i$, with
diagonal entries multiplied by a constant as presented in Table
\ref{dinwoodie}. The row sums and column sums of the entries are given
there too. Using {\tt LattE} we obtained the exact answer {\em
{8813835312287964978894}}.

\begin{table}
\begin{center}
\begin{tabular}{|c|c|c|c|c|c|c||c|} \hline
  2{\color{blue}{$x_1$}}  &   $x_2$   &   $x_3$ &  $x_4$ &  $x_5$ & $x_6$ & $x_7$ & 205\\ \hline
     $x_2$   &  2{\color{blue}{$x_8$}}   &   $x_{9}$ &  $x_{10}$ &  $x_{11}$ & $x_{12}$ & $x_{13}$ & 600 \\ \hline
     $x_3$  & $x_9$   &  2{\color{blue}{$x_{14}$}}  &  $x_{15}$ &  $x_{16}$ & $x_{17}$ & $x_{18}$ & 61\\ \hline
     $x_4$   & $x_{10}$ & $x_{15}$      &  2{\color{blue}{$x_{19}$}} &  $x_{20}$ & $x_{21}$ & $x_{22}$ & 17 \\ \hline
     $x_5$  & $x_{11}$ & $x_{16}$     &$x_{20}$  & 2{\color{blue}{$x_{23}$}}  & $x_{24}$ & $x_{25}$ & 11 \\ \hline
     $x_6$  &$x_{12}$  & $x_{17}$    & $x_{21}$  & $x_{24}$  & 2{\color{blue}{$x_{26}$}} & $x_{27}$ & 152 \\ \hline
     $x_7$ &  $x_{13}$ & $x_{18}$    & $x_{22}$  & $x_{25}$  & $x_{27}$  & 2{\color{blue}{$x_{28}$}} & 36 \\ \hline
    205 & 600 & 61 & 17 & 11 & 152 &36 & 1082 \\ \hline
\end{tabular}
\caption{The conditions for retinoblastoma RB1-VNTR genotype data from
 the Ceph database.} \label{dinwoodie} 
\end{center}
\end{table}

\end{document}